\documentclass[12pt]{amsart}
\usepackage{amsfonts}
\usepackage{amssymb}

\usepackage{amsmath}
\usepackage{epsf}
\usepackage{epsfig}
\usepackage{colordvi}
\usepackage{fancybox}

\usepackage{colordvi}

\newtheorem{thm}{Theorem}

\newcommand{\ds}{\displaystyle}

\begin{document}

\title[Characterizing a  surface by invariants]
{Characterizing a  surface by invariants}%


\author{ Ognian Kassabov}%

\address{Ognian Kassabov: University of Transport, Sofia, Bulgaria}
\email{okassabov@abv.bg}


\begin{abstract}
Canonical principal parameters are introduced for surfaces in  $\mathbb R^3$ without umbilical points. It is proved that
in these parameters the surface is determined (up to position in space) by a pair of invariants satisfying a
partial differential equation equivalent to the Gauss equation.  As such a pair of invariants we
may use the principal curvatures or the Gauss and the mean curvature.

\end{abstract}

\maketitle
\section{Introduction}
\label{intro}

An important problem in     differential geometry is to characterize a geometric object by its
invariants. For example, it is well known that  any curve in $\mathbb R^3$ is determined
(up to position in  space) by its curvature and torsion as functions of its natural parameter.

For the surfaces in $\mathbb R^3$ the situation is more complicated. According to the classical Bonnet's theorem,
a surface is determined (up to position in  space) by six functions -- the coefficients of the first and the second fundamental
forms satisfying the equations of Gauss and Codazzi. Of course  the coefficients of the fundamental forms are not invariant functions, 
unlike the curvature and torson of a curve, although these coefficients rest unchanged in motions. Nevertheless,
the above Bonnet's theorem helps us in studying the determination of a surface by invariants. 
Note that some differential equations between the invariants of the surfaces arise in a natural way as a result of
the equations of Gauss and Codazzi. The so-called Lund-Regge problem here is to find the minimum possible invariants 
and relations between them that characterize a surface, see  \cite{L-R}, \cite{Sym}. When
trying to reduce the number of invariants and the compatibility conditions, it is  common to
search for special parameters, just as in the case of the curves and their natural parameters.

An important progress in this direction is made  in \cite{G-M} -- a work 
that actually inspired the present paper. Namely in \cite{G-M}  it is  proved that a regular surface 
is determined (up to position in space) by four invariants -- the principal curvatures
$\nu_1$, $\nu_2$ and the geodesic curvatures $\gamma_1$, $\gamma_2$ of the principal lines. These invariants 
satisfy three partial differential equations equivalent to the Gauss and Codazzi equations.
In particular, for the class of Weingarten surfaces the authors  introduce some special
parameters that they call  {\it geometric} and they prove that in these parameters the surface is
determined by only one invariant function and two other functions. These three functions
are closely related to the principal curvatures and are subjects to a single
partial differential equation equivalent to the Gauss equation. 

In this paper, we introduce canonical parameters for any surface in $\mathbb R^3$ without umbilical points and
we prove that in these parameters the surface is locally determined up to position in  space by just two invariant functions 
 related by just one partial differential equation equivalent to the Gauss equation. 
These two invariant functions are the principal curvatures or the Gauss curvature and the mean curvature. It is clear that the surface cannot  
be determined by just one of these invariant functions -- for example there exist
many surfaces with the same constant  Gauss or mean curvature. 
So it appears that our results solve the Lund-Regge  problem for surfaces without umbilical points.

For similar investigations about surfaces in some upper dimensional  spaces of constant curvature $c$, we mention  
\cite{T-G}, where some special isothermal parameters are used in the case 
of minimal non-superconformal  surfaces in $\mathbb Q_c^4$ and it is proved that the surface is determined by the
Gauss curvature and the normal curvature, which satisfy a system of two partial differential equations;
see also \cite{G-K}.

\setcounter{equation}{0}
\section{Preliminaries}
\label{sec:1}

Let a regular surface in $\mathbb R^3$ be given by the parametric equation
$\	S \ :  \ x=x(u,v)  .\  $ 
We denote by $E$, $F$, $G$, resp. $L$, $M$, $N$
the coefficients of the first, resp. the second fundamental form. A point of $S$ is called 
{\it umbilical} if the two fundamental forms are proportional at that point. The Gauss curvature $K$ 
and the mean curvature $H$ of $S$, which are the most important invariants of the surface, 
are expressed with  these coefficients respectively by
$$
	K=\frac{LN-M^2}{EG-F^2} \qquad\qquad  H=\frac{EN-2FM+GL}{2(EG-F^2)} \ .
$$
Moreover the coefficients of the two fundamental forms 
satisfy the equation of Gauss
$$
	K=
	  -\frac1{2W}\left\{ \left(\frac{E_v-F_u}W\right)_v + \left(\frac{G_u-F_v}W\right)_u \right\}
	  -\frac1{4W^4}
	  \left|\begin{array}{ccc}
	    E   & F   & G \\
	    E_u & F_u & G_u \\
	    E_v & F_v & G_v
	    \end{array}\right|                   
$$
and the equations of Codazzi
$$
	2W^2(L_v-M_u)=(EN-2FM+GL)(E_v-F_u) + 
	\left|\begin{array}{ccc}
	    E   & F   & G \\
	    L   & M   & N \\
	    E_u & F_u & G_u 
	    \end{array}\right|
	 \ ,
$$
$$
	2W^2(M_v-N_u)=(EN-2FM+GL)(F_v-G_u) + 
	\left|\begin{array}{ccc}
	    E   & F   & G \\
			L   & M   & N \\
	    E_v & F_v & G_v 
	    \end{array}\right|
	 \ ,
$$
where $W=\sqrt{EG-F^2}$.
The classical theorem of Bonnet \cite{Bonnet} states that conversely, given six functions $E$, $F$, $G$, $L$, $M$, $N$
($E>0$, $EG-F^2>0$) that satisfy these equations, then locally there exists a unique (up to position in space) surface,
having  $E$, $F$, $G$ as coefficients of the first fundamental form and $L$, $M$, $N$ as
coefficients of the second fundamental form;  see also e.g. \cite{doCarmo}, p. 236.

Suppose a curve $c$ on $S$ be defined by
$$
	c\ :\ \ u=u(s)\ ,\ \ v=v(s) \ ,
$$ 
where $s$ is the natural parameter of $c$. Then the Frenet formulas are
$$
	\begin{array}{l}
		t'=\gamma p+\nu l \\
		p'=-\gamma t+\alpha l \\
		l'=-\nu t-\alpha p 
	\end{array}
$$  
where $t$ is the unit tangent vector field of $c$, $l$ is the unit normal vector field of $S$
and $p=l\times t$. The functions $\gamma$, $\nu$, $\alpha$ are respectively the geodesic curvature,
the normal curvature and the geodesic torsion of $c$ on $S$, respectively.
The normal curvature of $c$ is given by
$$
	\nu=\frac{L\dot u^2+2M\dot u\,\dot v+N\dot v^2}{E\dot u^2+2F\dot u\,\dot v+G\dot v^2}  \ .
$$
Actually at each point of $c$ the normal curvature $\nu$ depends not on the curve itself, but on the
direction of its tangent vector at that point, so we can speak about normal curvature of a direction in any point. 
The maximal and the minimal values of the normal curvatures at a point are called {\it principal curvatures}
and the corresponding directions and vectors -- {\it principal directions} and {\it principal vectors}. A curve on $S$
is called {\it principal} if its tangent vector is principal at any point.  
When the surface has no umbilical points, the parameters $(u,v)$ can be chosen such that
the parametric lines are principal. Then the parameters $(u,v)$ of $S$ are called {\it principal}. 
In terms of the coefficients of the fundamental forms
this means that  $F=M=0$ on $S$. 
In this case the geodesic torsions 
of the parametric lines vanish identically. On the other hand, the geodesic curvatures of
the parametric lines are
\begin{equation} \label{eq:2.1}
	\gamma_1=-\frac{E_v}{2E\sqrt G} 
	\ , \qquad
	\gamma_2=\frac{G_u}{2G\sqrt E} \ .
\end{equation}

Let $\nu_1$ and $\nu_2$ be the principal curvatures of $S$. Then the classical definition
of the Gauss curvature and the mean curvature becomes
\begin{equation} \label{eq:2.2}
	K=\nu_1\nu_2\ , \hspace {2cm} H=\frac12 (\nu_1+\nu_2) \ .
\end{equation}

\setcounter{equation}{0}
\section{Determining non-umbilical surfaces}
\label{sec:2}

Suppose that $S$ has no umbilical points and the
parametric lines are principal, i.e. $F=M=0$ on $S$. Then the  equation of Gauss is
\begin{equation} \label{eq:GaussEq}
	K=\nu_1\nu_2=
	  -\frac1{2\sqrt{EG}}\left\{ \left(\frac{E_v}{\sqrt{EG}}\right)_v +
											 \left(\frac{G_u}{\sqrt{EG}}\right)_u \right\}
\end{equation}
and the equations of Codazzi take the form
$$
	2EGL_v=(EN+GL)E_v \ ,    \qquad    2EGN_u=(EN+GL)G_u \ .
$$

On the other hand, the principal curvatures $\nu_1$, $\nu_2$  are given by
\begin{equation} \label{eq:3.2}
	\nu_1=\frac LE\ , \qquad \nu_2=\frac NG \ .
\end{equation}
Since the surface has no umbilical points, the difference  $\nu_1-\nu_2$ cannot vanish. Hence 
it is easy to see that the equations of Codazzi  may be written as
\begin{equation} \label{eq:CodazziEqForNu}
	\frac{E_v}{2E}=-\frac{(\nu_1)_v}{\nu_1-\nu_2} \ ,    \qquad \qquad       \frac{G_u}{2G}=\frac{(\nu_2)_u}{\nu_1-\nu_2} \ .
\end{equation}
Let us fix a point $(u_0,v_0)$. The last equations imply that  there exist two functions $\varphi_1(u)$ and $\varphi_2(v)$, such that
$$
	\sqrt E=\varphi_1(u)e^{-\ds\int_{v_0}^v\frac{(\nu_1)_v}{\nu_1-\nu_2}dv} \ ,    \qquad \qquad       
	\sqrt G=\varphi_2(v)e^{\ds\int_{u_0}^u\frac{(\nu_2)_u}{\nu_1-\nu_2}du} \ .
$$
In other words, for any functions $\phi_1(u)$, $\phi_2(v)$, 
the function
$$
	\phi_1(u)\sqrt E\,e^{\ds\int_{v_0}^v\frac{(\nu_1)_v}{\nu_1-\nu_2}dv}
$$
does not depend on $v$ and the function
$$
	\phi_2(v)\sqrt G\,e^{-\ds\int_{u_0}^u\frac{(\nu_2)_u}{\nu_1-\nu_2}du}
$$
does not depend on $u$.
Now we introduce new parameters $(\bar u,\bar v)$ by the formulas
$$
	\begin{array}{l}
		\ds\overline u=\frac1{\sqrt{E(u_0,v_0)}}\int_{u_0}^u \sqrt E\, e^{\ds\int_{v_0}^v\frac{(\nu_1)_v}{\nu_1-\nu_2}dv +
		                                       \int_{u_0}^u\frac{(\nu_1)_u}{\nu_1-\nu_2}(u,v_0)du } du + \overline u_0\ , \\
		\ds\overline v=\frac1{\sqrt{G(u_0,v_0)}}\int_{v_0}^v \sqrt G\,e^{-\ds\int_{u_0}^u\frac{(\nu_2)_u}{\nu_1-\nu_2}du -
		                                       \int_{v_0}^v\frac{(\nu_2)_v}{\nu_1-\nu_2}(u_0,v)dv}dv  +\overline v_0\ .
	\end{array}
$$
for some constants $\overline u_0$, $\overline v_0$. The parameters $(\overline u,\overline v)$ are also principal.  Moreover we have 
\begin{equation} \label{eq:canon}
	\begin{array}{l}
		\ds\frac{\sqrt{\overline E}}{\sqrt{\overline E(\overline u_0,\overline v_0)}}\,e^{\ds 
		\int_{\overline v_0}^{\overline v}\frac{(\overline \nu_1)_{\overline v}}{\overline \nu_1-\overline \nu_2}d\overline v+
		\int_{\overline u_0}^{\overline u}\frac{(\overline \nu_1)_{\overline u}}{\overline \nu_1-\overline \nu_2}(\overline u,\overline v_0)d\overline u}=1     \\     
		\ds\frac{\sqrt{\overline G}}{\sqrt{\overline G(\overline u_0,\overline v_0)}}\,e^{\ds 
		-\int_{\overline u_0}^{\overline u}\frac{(\overline \nu_2)_{\overline u}}{\overline \nu_1-\overline \nu_2}d\overline u-
		\int_{\overline v_0}^{\overline v}\frac{(\overline \nu_2)_{\overline v}}{\overline \nu_1-\overline \nu_2}(\overline u_0,\overline v)d\overline v} =1
	\end{array} \ .
\end{equation}
We shall call {\it canonical principal parameters} any principal parameters $(\overline u,\overline v)$ satisfying  (\ref{eq:canon})
for certain constants $(\overline u_0, \overline v_0)$.         
 
We can see by a straightforward check that if $(u,v)$ are also canonical principal parameters,
then
$$
	\begin{array}{l}
	 \overline  u=\lambda\, u +c_1 \\
	 \overline  v=\mu\, v+c_2 
	\end{array}
	\qquad\qquad {\rm or} \qquad\qquad
	\begin{array}{l}
	 \overline  u=\lambda\, v +c_1 \\
	 \overline  v=\mu\, u+c_2 
	\end{array}
$$
for some constants $\lambda$, $\mu$, $c_1$, $c_2$ ($\lambda\ne 0$, $\mu\ne 0$). More precisely, if 
$\overline  u_0=\overline u(u_0)$,  $\overline  v_0=\overline v(v_0)$, then
$$
	\begin{array}{l}
	(\overline  u-\overline  u_0) \sqrt{\overline   E(\overline  u_0,\overline  v_0)}=(u-u_0) \sqrt{E(u_0,v_0)} \\
	(\overline  v-\overline  v_0) \sqrt{\overline   G(\overline  u_0,\overline  v_0)}=(v-v_0) \sqrt{G(u_0,v_0)} 
	\end{array}\ .
$$
        
In the following we assume that the surface is parametrized with canonical  principal para\-meters $(u,v)$. Then
the coefficients $E$ and $G$ of the first fundamental form satisfy

\begin{equation} \label{eq:3.3}
	\begin{array}{l}
		E=a\,e^{-\ds 2\int_{v_0}^v\frac{(\nu_1)_v}{\nu_1-\nu_2}dv-2\int_{u_0}^u\frac{(\nu_1)_u}{\nu_1-\nu_2}(u,v_0)du}     \\     
		G=b\,e^{\ds 2\int_{u_0}^u\frac{(\nu_2)_u}{\nu_1-\nu_2}du+2\int_{v_0}^v\frac{(\nu_2)_v}{\nu_1-\nu_2}(u_0,v)dv} 
	\end{array} \ ,
\end{equation}
where $a=E(u_0,v_0)$, $b=G(u_0,v_0)$. In this case the Gauss equation (\ref{eq:GaussEq}) can be written in the following equivalent form

\begin{equation} \label{eq:NewGaussEq}
	\nu_1 \nu_2\Psi_1\Psi_2=\frac1b\left( \frac{(\nu_1)_v}{\nu_1-\nu_2}\, \frac{\Psi_1}{\Psi_2}\right)_v 
											-   \frac1a\left(\frac{(\nu_2)_u}{\nu_1-\nu_2}\,  \frac{\Psi_2}{\Psi_1}\right)_u 	                             
\end{equation} 
where the functions $\Psi_1$ and $\Psi_2$ are defined by
\begin{equation} \label{eq:Psi1Psi2}
	\begin{array}{l}
		\Psi_1=\,e^{-\ds \int_{v_0}^v\frac{(\nu_1)_v}{\nu_1-\nu_2}dv-\int_{u_0}^u\frac{(\nu_1)_u}{\nu_1-\nu_2}(u,v_0)du}     \\     
		\Psi_2=\,e^{\ds \int_{u_0}^u\frac{(\nu_2)_u}{\nu_1-\nu_2}du+\int_{v_0}^v\frac{(\nu_2)_v}{\nu_1-\nu_2}(u_0,v)dv} 
	\end{array} \ .
\end{equation}

Conversely, consider two differentiable functions $\nu_1$, $\nu_2$ that satisfy the
equation (\ref{eq:NewGaussEq}) for some positive constants $a$, $b$, the functions $\Psi_i$  being defined by   (\ref{eq:Psi1Psi2})
(of course we suppose that the difference $\nu_1-\nu_2$ never vanishes). 
With these functions $\nu_1$, $\nu_2$ we define $E$ and $G$ by
(\ref{eq:3.3}) and after that  $L$ and $N$ by (\ref{eq:3.2}).
Then using the theorem of Bonnet we obtain:

\begin{thm}\label{T:MainTheorem1}  
Let $a$ and $b$ be positive constants and two differentiable functions $\nu_1(u,v) $, $\nu_2 (u,v) $ be given.
Define $\Psi_1$, $\Psi_2$ by (\ref{eq:Psi1Psi2}) and suppose that 
(\ref{eq:NewGaussEq}) is satisfied. Then locally there exists a unique 
(up to position in  space) surface $S$, such that   $\nu_1$ and $\nu_2$
are the principal curvatures of $S$ in canonical principal parameters. For this
surface $E(u_0,v_0)=a$, $G(u_0,v_0)=b$.
\end{thm}

Note  that  the integrability condition (\ref{eq:NewGaussEq}) (which is a form of the
Gauss equation) is expressed only by the two invariants $\nu_1$ and $\nu_2$ -- the
principal curvature functions of the surface in canonical principal parameters.

Note also that the above theorem and the Gauss integrability equation (\ref{eq:NewGaussEq})
can be put in a different form in terms of the Gauss curvature and the normal curvature
instead of the principal curvatures $\nu_1$, $\nu_2$. Indeed according to (\ref{eq:2.2})
we have (up to numeration)
$$
	\nu_1=H+\sqrt{H^2-K}\ ,  \hspace{2cm} \nu_2=H-\sqrt{H^2-K}\ .
$$
In this case the condition that $\nu_1-\nu_2$ never vanishes is replaced by the
condition that $H^2-K$ never vanishes. As a result, the surface is
determined up to position in space by its Gauss and mean curvature.
More precisely, we obtain

\begin{thm}\label{T:MainTheorem2}  
Let $K(u,v)$, $H(u,v)$ be differentiable functions such that the 
equation 
$$
		\ds\frac{2K}{\sqrt{H^2-K}}\, \Phi_1\Phi_2 
			=\ds	 \frac 1b	\left(\frac{\Phi_1}{\Phi_2}\,\frac{\Big(H+\sqrt{H^2-K}\Big)_v}{\sqrt{H^2-K}}  \right)_v
			- \ds  \frac 1a\left(\frac{\Phi_2}{\Phi_1}\,\frac{\Big(H-\sqrt{H^2-K}\Big)_u}{\sqrt{H^2-K}}  \right)_u 
$$
where
$$
	\begin{array}{l}
		\Phi_1=\,e^{-\ds \int_{v_0}^v\frac{H_v}{2\sqrt{H^2-K}}dv- \int_{u_0}^u\frac{H_u}{2\sqrt{H^2-K}}(u,v_0)du}     \\     
		\Phi_2=\,e^{\ds  \int_{u_0}^u\frac{H_u}{2\sqrt{H^2-K}}du+\int_{v_0}^v\frac{H_v}{2\sqrt{H^2-K}}(u_0,v)dv} 
	\end{array} 
$$
is satisfied for some positive constants $a$, $b$. Then locally there exists a unique 
(up to position in  space) surface, such that   $K$ and $H$ are respectively 
its Gauss curvature and mean curvature in canonical principal parameters. 
For this surface $\big(E\sqrt{(H^2-K)}\,\big)(u_0,v_0)=a$, $\big(G\sqrt{(H^2-K)}\,\big)(u_0,v_0)=b$.
\end{thm}

Having two functions $\nu_1$, $\nu_2$ satisfying the conditions of Theorem \ref{T:MainTheorem1} (or, what is
the same, two functions $K$, $H$ satisfying the conditions of Theorem \ref{T:MainTheorem2}),
we determine the coefficients $E$, $G$ of the first fundamental form of the induced surface $S$ by  (\ref{eq:3.3}). Now
we can find the geodesic curvatures $\gamma_1$, $\gamma_2$ of the principal lines of the surface using (\ref{eq:2.1}).
A geometric method to construct the surface with invariants  $\nu_1$, $\nu_2$, $\gamma_1$, $\gamma_2$
is obtained in \cite{G-M}.

\setcounter{equation}{0}
\section{Particular cases}
\label{sec:3}

The surface $S\ :\ x=x(u,v),\ (u,v)\in D$ is called {\it strongly regular Weingarten surface} if
$$
	\big(\nu_1(u,v)-\nu_2(u,v)\big) \gamma_1(u,v)\gamma_2(u,v) \ne 0, \qquad (u,v)\in D
$$                                                                                 
and there exist two differentiable functions $f(t)$, $g(t)$ defined on an interval $I$ 
and a function $\nu(u,v)$, defined on $D$, such that
\begin{equation} \label{eq:4.1}
	f(t)-g(t)>0\,, \qquad f'(t)g'(t) \ne 0 \,, \ \ t\in I\,,
\end{equation}
\begin{equation} \label{eq:4.2}
	 \nu_u(u,v)\nu_v(u,v)\ne 0, \ \ (u,v)\in D \ .
\end{equation}	
\begin{equation} \label{eq:4.3}
	\nu_1=f(\nu)\,, \quad \nu_2=g(\nu)\,.
\end{equation}	
Theorem \ref{T:MainTheorem1} implies that given  three functions $f(t)$, $g(t)$, $\nu(u,v)$ with the  
properties (\ref{eq:4.1}), (\ref{eq:4.2}) and satisfying the equation
\begin{equation} \label{eq:4.4}
	\begin{array}{rcl}
		\ds A\, \left\{ f'\nu_{vv}+\left(f''-\frac{2f'^2}{f-g}\right) \nu_v^2 \right\}
	               e^{\ds 2\int_{\nu_0}^{\nu} \frac{g'dt}{g-f} } &  & \\
		-\ds  B\,  \left\{ g'\nu_{uu}+\left(g''+\frac{2g'^2}{f-g}\right) \nu_u^2 \right\}
	              e^{\ds  2\int_{\nu_0}^{\nu} \frac{f'dt}{f-g} } & = &  fg(f-g) 							
	\end{array} 
\end{equation}
for two positive constants $A$, $B$ and $\nu_0=\nu(u_0,v_0)$ for $(u_0,v_0) \in D$, 
then there exists a unique (up to position in space) 
Weingarten surface $S$ with principal curvatures in canonical principal parameters given by (\ref{eq:4.3}). 
This is one of the main results in \cite{G-M}. Note that in this case our canonical principal parameters
coincide with the geometric  principal parameters, defined in  \cite{G-M}.

For the form of the Gauss equation (\ref{eq:4.4}) for some 
important subclasses of Wiengarten surfaces, e.g. surfaces of constant mean curvature, see \cite{G-M}.

It is more interesting to consider the surfaces of constant mean curvature $H$ by another point of view.
Namely, according to Theorem \ref{T:MainTheorem2} such a surface is uniquely 
determined  by its Gauss curvature. More precisely  Theorem \ref{T:MainTheorem2} (with $a=b=1$) implies that for a real number $H$ and a
differentiable function $K$ satisfying $K<H^2$ and the differential equation
$$
	\ds\Delta(\log(H^2-K))=	\frac{4K}{\sqrt{H^2-K}}
		 \,,
$$
where $\Delta$ is  the Laplace operator, there exists a unique (up to position) surface with 
Gauss curvature $K$ and constant mean curvature $H$.

In particular, for minimal surfaces $(H=0)$ this equation reduces to
$$
	\ds\Delta\big(\log\sqrt{-K}\big)+2	\sqrt{-K}=0
$$
or, if $\nu=\sqrt{-K}$ is the positive principal curvature,
\begin{equation} \label{eq:4.5}
	\Delta(\log\nu)+	2\nu=0 \,.
\end{equation}
According to (\ref{eq:3.3}), in this case $E=G$ and since $F=0$, the canonical principal parameters $(u,v)$ are isothermal.  
When we consider minimal surfaces, it is very common to use complex coordinates; in  
real coordinates this gives isothermal parameters. A method to obtain
canonical principal parameters for a minimal surface from arbitrary isothermal ones is found in \cite{OK}.
In \cite{Ganchev} the equation (\ref{eq:4.5}) is named {\it natural partial differential equation of minimal surfaces}. 

The flat surfaces, i.e. the surfaces with vanishing Gauss curvature $K$, are well studied -- they are
general cylinders, general cones and tangent developable surfaces.
When the surface has no umbilical points (for example for a tangent developable surface 
the torsion of the directrix must not vanish) the mean curvature $H $ can not vanish. It follows 
from  Theorem \ref{T:MainTheorem2} that these surfaces are characterized by
$$
	\left(\frac{1}{H}\right)_{vv}=0
	\qquad {\rm or} \qquad \ds H=\frac1{f(u)v+g(u)}
$$
in canonical principal parameters for some functions $f(u)$, $g(u)$.

\end{document}